# Dynamic robust stabilization of fractional-order linear systems with nonlinear uncertain parameters: An LMI approach

Pouya Badri[1], Mahdi Sojoodi[1], Elyar Zavvari[1]
[1]Advanced Control Systems Laboratory, School of Electrical and Computer Engineering, Tarbiat Modares University, Tehran, Iran.

**Abstract: This paper presents a dynamic output feedback controller with determined order for the stabilization of a class of fractional-order system with nonlinear uncertain parameters with fractional order $0 < \alpha < 2$. Using stability theories of fractional-order systems and linear matrix inequalities, some sufficient conditions in the LMI form are deduced to guarantee the robustness and asymptotic stabilization of the system. Designing a dynamic robust controller, along with all its useful features, leads to more unknown parameters in comparison with a static controller and makes controller design procedure more difficult due to more complex constraints taht must be solved. In this paper, using proper lemmas and theorems, LMI techniques, and suitable solvers and parsers the difficulty of designing such controllers has been overcome. Simulation results of three different numerical examples illustrate that the proposed sufficient theoretical results are applicable and effective for tackling robust stabilization problems.**

**Keywords: Fractional-order system, nonlinear uncertain parameters, linear matrix inequality (LMI), robust stabilization, dynamic output feedback.**

## 1. Introduction

In the last decades, fractional-order calculus has received considerable interest and attention of physicists and engineers and have found many applications in various fields such as viscoelastic, electrode-electrolyte, biological electric conductance, neural systems, and others (Badri and Sojoodi, 2019b; Chen et al., 2017; Badri and Saleh Tavazoei, 2016; Liang et al., 2019; Xie et al., 2019; Soukkou et al., 2018; Sumathi and Umasankar, 2018). Despite integer-order derivative which describes local properties of a certain position, or a variation at a specific moment in time for a physical process, fractional-order derivative is related to the whole space and the whole time domain (Badri and Tavazoei, 2016). Therefore, fractional-order differential equations can more completely and precisely describe systems having responses with long memory transients than the ordinary integer-order differential equations (Badri and Tavazoei, 2014, 2017). Accordingly, stability and stabilization of fractional-order systems is an important and challenging problem since many physical and real-world processes are modeled with fractional-order state equations (Badri and Sojoodi, 2018; Ivanova et al., 2018; Lu and Chen, 2009; Ma et al., 2014; Zhang et al., 2018; Badri and Sojoodi, 2019a).

Unfortunately, uncertainties arising from neglected dynamics, uncertain physical parameters, parametric variations in time, and many other sources are inevitable in real physical system. Thus, robust stability and stabilization problems have become an essential issue for all control systems including fractional-order ones (Badri and Sojoodi, 2018; Lu and Chen, 2009; Ma et al., 2014). Robust stability and stabilization problems of fractional-order systems were investigated in (Badri and Sojoodi, 2018; Lu and Chen, 2009; Ma et al., 2014; Xing and Lu, 2009; Chen and Lin, 2004; Alaviyan Shahri et al., 2018; Binazadeh and Yousefi, 2018; Dadras and Momeni, 2014; Boubellouta et al., 2019). Stability and stabilization problem of FO-LTI interval systems was investigated in (Lu and Chen, 2009), in which necessary and sufficient conditions were presented in LMI form. Furthermore, necessary and sufficient conditions for checking robust stability of general interval FO-LTI system were investigated in (Zheng, 2017), in which interval uncertainties exist both in the coefficients and orders of the system. In (Ma et al., 2014) the robust stability and stabilization problems of fractional-order linear systems with positive real uncertainty were solved, where the existence conditions and design procedures of the static state feedback, static output feedback and observer-based controllers for asymptotically stabilizing of such systems were investigated with the constraint on the output matrix to be of full-row rank. Moreover, the stability and stabilization problem for a class of uncertain fractional-order systems subject to input saturation was investigated in (Alaviyan Shahri et al., 2018). Besides, robust stabilization for a class of nonlinear time-delay fractional-order systems in the presence of nonlinear Lipschitz functions; time-varying norm-bounded uncertain terms; and time-delays in the state variables was studied in (Binazadeh and Yousefi, 2018).

Linear fractional-order systems have been used in modeling of a wide range of fractional-order systems such as linear electrical circuits with fractional active elements, civil structures, and etc. On the other hand, most of nonlinear fractional-order systems can be linearized in practical applications for stability analysis and control purposes. Therefore, the stability of fractional-order linear systems has received increasing attention and interest in recent years (Ma et al., 2014; Badri and Sojoodi, 2019b; Zhang and Zhao, 2020; Ghorbani, 2020). Motivated by above discussion, robust stabilization of fractional-order linear systems with nonlinear uncertainty has been investigated in this paper.

The normalization and stabilization of rectangular descriptor fractional-order interval systems with fractional order $0 < \alpha < 1$ is considered in (Zhang and Zhao, 2020), where a rectangular descriptor fractional order interval system is transformed into an augmented square descriptor fractional-order interval system by adopting the proportional and derivative type dynamic compensator. Furthermore, in (Zhang and Wang, 2020), the stability and robust stabilization of switched fractional-order systems are concerned. Firstly, two stability theorems for switched fractional-order systems with order $0<\alpha<1$ and $1<\alpha<2$ under the arbitrary switching law are given. Secondly, the relationship between the stability of switched integer order systems and that of switched fractional-order systems is obtained. Moreover, robust stability of the FO-LTI systems with poly-topic and two-norm bounded uncertainties of the fractional order $1 < \alpha < 2$ was investigated in (Li, 2018) and the state feedback controller was designed for such systems as well. Besides, in (Ghorbani, 2020) the robust stability of fractional-order plants suffering from interval uncertainties by using fractional-order controllers is investigated based on the zero-exclusion principle, in which, new necessary and sufficient criteria are proposed to analyze the robust stability of the corresponding characteristic polynomial.

The robust stability of fractional-order systems described in the pseudo-state space model with incommensurate fractional orders is investigated in (Tavazoei and Asemani, 2020), in which a stability criterion for integer-order systems is extended to incommensurate-order fractional systems by using the generalized Nyquist Theorem.

In most of the mentioned works state feedback controller is employed. In state feedback control scheme, all individual states of the system are needed to be measured and used in feedback line. However, in some cases measuring all states is impossible due to economic issues or physical limitations, where using output feedback control seems to be effective. In output feedback scheme, there is no need to measure all individual states of the system and only by measuring outputs of the system the control action is done (Badri et al., 2016, 2019).

Besides, it has been claimed that dynamic feedback controller brings about more effective control performances, flexibility, and degrees of freedom for the sake of achieving control objectives, in comparison with the static one (Park, 2009). Furthermore, it can be easily shown that some unstable systems cannot be stabilized by static controllers and a dynamic controller is merely needed to ensure the stabilization of such systems (Sontag, 1998). Then, by means of linear matrix inequalities, we have designed a robust dynamic output feedback controller for FO-LTI systems with positive real uncertainty in (Badri and Sojoodi, 2018).

In (Lan and Zhou, 2013), robust stabilization problem of a class of Lipschitz non-linear fractional-order systems with fractional order $0 < \alpha < 1$, was investigated using an observer-based robust stabilization controller, where the stabilization results are reliable if the output matrix of the uncertain system is of full row rank. Robust stability and stabilization sufficient conditions were derived in (Ibrir and Bettayeb, 2015) to stabilize fractional-order systems, with fractional order $0 < \alpha < 1$, subject to bounded uncertainties using state feedback and observer-based controllers. In this paper, it is assumed that all individual possible pairs of $A + \Delta A$, $C + \Delta C$ are observable in the sense of Kalman while the system uncertainties are randomly distributed in the state matrix $A$ and the output matrix $C$. Robust D-stability test of LTI general fractional order control systems in a closed-loop is investigated in (Mohsenipour and Liu, 2020), in which either of the system or the controller may be of fractional order and a necessary and sufficient condition for testing the robust D-stability of these systems is derived.

It is worth noting that the abovementioned uncertainty models allow only linearly dependent uncertain parameters. Nevertheless, practically, in a real system, the uncertain parameters usually appear in a nonlinear form (e.g. geometric and inertia parameters for a dynamic system) (Xu and Darouach, 1998). The above uncertainty descriptions may seriously result in ''over bounds" of uncertainties, which may lead to extremely conservative robustness analysis results (Xu and Darouach, 1998). Therefore, nonlinear uncertainty has attracted the attention of many researchers in the last years (Amini et al., 2016; Badri et al., 2016; Jianbin et al., 2020).

In most of available controller design methods, high-order controllers are produced suffering from costly implementation, high fragility, unfavorable reliability, maintenance problems, and possible numerical errors. Since the desired closed-loop performance is not always assured by using available plant or controller order reduction methods, designing a low-order controller would be helpful in this case, i.e., the order of the controller is predetermined in advance, (Badri et al., 2016). To the best knowledge of authors, there is no research on the analytical design of a stabilizing fixed-order dynamic output feedback controller for fractional-order systems with nonlinear uncertain parameters.

Motivated by mentioned observations, the purpose of our paper is to solve the problem of robust stabilization of fractional-order linear system with nonlinear uncertain parameters with the fractional order $0 < \alpha < 2$, using a fixed-order dynamic output feedback controller, where sufficient conditions are presented by means of linear matrix inequalities (LMIs). Using a predetermined order dynamic output feedback controller makes it possible to use a low-order controller to stabilize the uncertain system. Furthermore, in our proposed method designing a dynamic feedback controller does not bring about limiting constraints on the state-space matrices of the uncertain systems which is usual in most of the previous methods. Despite the complexity of considering the most complete model of linear controller containing direct feedthrough parameter, the LMI form of the constraints is preserved, making them suitable to be used in practice thanks to various efficient convex optimization parsers and solvers that can be applied to determine the feasibility of constraints and calculate design parameters.

To the best of our knowledge, there is no result on the analytical design of a stabilizing fixed-order dynamic output feedback controller for fractional-order systems with nonlinear uncertain parameters in the literature. Utilizing dynamic controller, more performance efficiency, and degrees of freedom can be achieved. Besides, not all of the pseudo-states of the uncertain FO-LTI system is necessary, thanks to the output feedback scheme. In this paper, sufficient conditions are presented for designing a robust stabilizing controller with a predetermined order, which can be chosen to be as low as possible for simpler implementation and lower cost. Furthermore, A unified method is proposed to stabilize FO-LTI systems with the non-integer order $\alpha$ satisfying $0 < \alpha < 1$ or $1 \leq \alpha < 2$.

The structure of the paper is as follows. In section 2, some necessary preliminaries and lemmas together with the problem formulation are presented. Robust stabilizing conditions of fractional-order uncertain systems via a dynamic output feedback controller are derived in Section 3. Section 4 presents some numerical examples to verify the results. Finally, conclusions are drawn in section 5.

**Notations**: In this paper, $A \otimes B$ denotes the Kronecker product of two matrices, and $M^T$, $\overline{M}$, and $M^*$ stand respectively, for the transpose, the conjugate, and the transpose conjugate of $M$. The conjugate of the scalar number $z$ is represented by $\bar{z}$ and $Sym(M)$ is used to denote the expression $M + M^*$. The notation $*$ is the symmetric component symbol in matrix and $\uparrow$ is the symbol of pseudo-inverse of matrix. The notations **0** and $I$ denote the zero and identity matrices with appropriate dimensions and $i$ is used to represent the imaginary unit.

## 2. Preliminaries and problem formulation

In this paper, the following Caputo definition for α fractional order derivatives of function $f(t)$ is adopted since the Laplace transform of this definition copes with clear physical interpretations [28]:

$$ {}_{a}^{C}D_{t}^{\alpha}f(t) = \frac{1}{\Gamma(m-\alpha)}\int_{a}^{t}(t-\tau)^{m-a-1}\left(\frac{d}{d\tau}\right)^{m} f(\tau)d\tau $$

where $\Gamma(\cdot)$ is Gamma function defined by $\Gamma(\epsilon) = \int_{0}^{\infty} e^{-t}t^{\epsilon-1}dt$ and $m$ is the smallest integer that is equal to or greater than $\alpha$.

Consider the following generalized uncertain fractional-order mathematical model of the integer-order one in (Xu and Darouach, 1998), which is more prevalent in real physical systems in comparison with models containing linearly dependent uncertain parameters.

$$ (A_n + \Delta_n)D^{n\alpha}q + (A_{n-1} + \Delta_{n-1})D^{(n-1)\alpha}q + \cdots + (A_1 + \Delta_1)D^{\alpha}q + (A_0 + \Delta_0)q = E \tag{1} $$

where $\alpha$ is fractional order, $A_i \in R^{m\times m}$ and $\Delta_i \in R^{m\times m}$ are respectively known matrices which represent the values of the system at the nominal working point and the unknown matrix representing the uncertain parameters. Vector $E$ represents a known driving source. The $i\alpha$th order differential of vector $q$ is represented by $D^{i\alpha}q,\ i = 1,2,\dots,n$. Furthermore, as indicated in (Xu and Darouach, 1998), $A_n + \Delta_n$ is assumed to be nonsingular which is true for most of the physical systems. Equation (1) can be rewritten as the following compact form:

$$ \begin{aligned} &D^{\alpha}x(t) = (I + \Delta I)[(A + \Delta A)x(t) + Bu(t)], 0 < \alpha < 2 \\ &y(t) = Cx(t) \end{aligned} \tag{2} $$

where

$$ x = \begin{bmatrix} q \\ D^{\alpha}q \\ \vdots \\ D^{(n-1)\alpha}q \end{bmatrix}, \Delta_I = \begin{bmatrix} 0 & \cdots & 0 & 0 \\ \vdots & \ddots & \vdots & \vdots \\ 0 & \cdots & 0 & 0 \\ 0 & \cdots & 0 & -(A_n + \Delta_n)^{-1}\Delta_n \end{bmatrix}, A = \begin{bmatrix} 0 & I_m & \cdots & 0 \\ \vdots & \vdots & \ddots & \vdots \\ 0 & 0 & \cdots & I_m \\ -A_n^{-1}A_0 & -A_n^{-1}A_1 & \cdots & -A_n^{-1}A_{n-1} \end{bmatrix}, $$

$$ \Delta_A = \begin{bmatrix} 0 & 0 & \cdots & 0 \\ \vdots & \vdots & \ddots & \vdots \\ 0 & 0 & \cdots & 0 \\ -A_n^{-1}\Delta_0 & -A_n^{-1}\Delta_1 & \cdots & -A_n^{-1}\Delta_{n-1} \end{bmatrix},\ B = \begin{bmatrix} 0 \\ \vdots \\ 0 \\ A_n^{-1} \end{bmatrix}, u(t) = E, \tag{3} $$

in which $x \in R^n$ denotes the pseudo-state vector, $u \in R^l$ is the control input, $y \in R^m$ is the output vector, and $C \in R^{m\times n}$ is the output matrix. Furthermore, uncertainties $\Delta_I$ and $\Delta_A$ satisfy the following polytopic structure

$$ \Delta_I = \sum_{i=1}^{p} m_i M_i, \Delta_A = \sum_{i=1}^{q} n_i N_i, \tag{4} $$

where $m_i$ and $n_i$ are the uncertain parameters satisfying the bounds $|m_i| \le m$, and $|n_i| \le n$, with constants $m, n > 0$. Also, $M_i$ and $N_i$ are known constant matrices with proper dimensions. It can be easily obtained from (4) that

$$ \Delta_I\Delta_I^T \le H,\ \Delta_A\Delta_A^T \le G, \tag{5} $$

with $H = \sum_{i=1}^{p} pm^2 M_i M_i^T$ and $G = \sum_{i=1}^{q} qn^2 N_i N_i^T$.

The following lemmas are needed in order to study the stability of fractional-order systems and obtain the main results.

**Lemma 1** (Farges et al., 2010): Let $A \in R^{n\times n}$, $0 < \alpha < 1$ and $\theta = (1-\alpha)\pi/2$. The fractional-order system $D^{\alpha}x(t) = Ax(t)$ is asymptotically stable if and only if there exists a positive definite Hermitian matrix $X = X^* > 0, X \in C^{n\times n}$ such that

$$ (rX + \bar{r}\bar{X})^T A^T + A(rX + \bar{r}\bar{X}) < 0, \tag{6} $$

where $r = e^{\theta i}$.

**Lemma 2** (Sabatier et al., 2010): Let $A \in \mathcal{R}^{n\times n}$, $1 \le \alpha < 2$ and $\theta = \pi - \alpha\pi/2$. The fractional-order system $D^{\alpha}x(t) = Ax(t)$ is asymptotically stable if and only if there exists a positive definite matrix $X \in \mathcal{R}^{n\times n}$ such that

$$ \begin{bmatrix} (A^TX + XA)sin\,\theta & (XA - A^TX)cos\theta \\ * & (A^TX + XA)sin\,\theta \end{bmatrix} < 0, \tag{7} $$

defining

$$ \Theta = \begin{bmatrix} \sin\theta & -\cos\theta \\ \cos\theta & \sin\theta \end{bmatrix}, \tag{8} $$

and with this in mind that $A$ is similar to $A^T$, inequality (15) can be stated as follows

$$ Sym\{\Theta \otimes (AX)\} < 0. \tag{9} $$

**Lemma 3** (Lu and Chen, 2009): For any matrices $X$ and $Y$ with appropriate dimensions, we have

$$X^TY + Y^TX \leq \eta X^TX + (1/\eta)Y^T\,Y\ for\ any\ \eta > 0. \tag{10}$$

**Lemma 4** (Li and De Souza, 1997): For any real matrices $X \in R^{n\times n}$, $Y \in R^{n\times n}$, and scalar $\varepsilon > 0$ such that $I - \varepsilon YY^T > 0$, we have

$$(X + Y)^T(X + Y) \leq X^T(I - \varepsilon YY^T)^{-1}X + \varepsilon^{-1}I. \tag{11}$$

## 3. Main results

In this section, first, a new display for the given FO-LTI system (2) is given. Then, a new stabilization condition is derived for the uncertain system and an LMI approach is proposed for designing a dynamic output feedback control law to stabilize the fractional-order system (2), using the stabilizing controller parameters.

The uncertain FO-LTI system (2) can be rewritten as follows

$$\begin{aligned} &D^\alpha x(t) = \tilde{A}x(t) + \tilde{B}u(t), 0 < \alpha < 2 \\ &y(t) = Cx(t) \end{aligned} \tag{12}$$

with

$$\tilde{A} = (I + \Delta_I)(A + \Delta_A), \tilde{B} = (I + \Delta_I)B. \tag{13}$$

The main aim of this paper is to design a robust dynamic output feedback controller that asymptotically stabilizes the FO-LTI system (2) with nonlinear uncertain parameters $\Delta_I$ and $\Delta_A$ in the terms of LMIs. Therefore, the following dynamic output feedback controller is presented

$$\begin{aligned} &D^\alpha x_C(t) = A_Cx_C(t) + B_Cy(t), \qquad 0 < \alpha < 2 \\ &u(t) = C_Cx_C(t) + D_Cy(t), \end{aligned} \tag{14}$$

with $x_C \in \mathcal{R}^{n_c}$, in which $n_c$ is the arbitrary order of the controller and $A_C, B_C, C_C$, and $D_C$ are appropriate matrices to be designed. The resulted closed-loop augmented system employing (2) and (14) is as follows

$$D^\alpha x_{Cl}(t) = \tilde{A}_{Cl}x_{Cl}(t), \qquad 0 < \alpha < 2 \tag{15}$$

where

$$x_{Cl}(t) = \begin{bmatrix} x(t) \\ x_C(t) \end{bmatrix}, \quad \tilde{A}_{Cl} = \begin{bmatrix} \tilde{A} + \tilde{B}D_CC & \tilde{B}C_C \\ B_CC & A_C \end{bmatrix}. \tag{16}$$

In the following, two robust stabilizing theorems for uncertain FO-LTI system (12) are given in terms of LMIs for $0 < \alpha < 1$ and $1 \leq \alpha < 2$ cases, respectively.

**Theorem 1**: Considering closed-loop system in (15) with $0 < \alpha < 1$, and a positive definite Hermitian matrix $P = P^*$ in the form of

$$P = diag(P_S, P_C), \tag{17}$$

with $P_S \in C^{n\times n}$ and $P_C \in C^{n_c\times n_c}$ and real positive scalar constants $\eta_1$, $\varepsilon_1$, and $\eta_3$ alongside with matrices $T_i, i = 1, \dots, 4$ exist such that the following LMI constrain become feasible

$$\begin{bmatrix} \Sigma_{11} + \eta_1\tilde{G} + (\eta_3 + 1)\tilde{H} & (rP + \bar{r}\bar{P})^T & \Sigma_{13} & \Sigma_{14} \\ * & -\eta_1 I & \mathbf{0} & \mathbf{0} \\ * & * & \Sigma_{33} & \mathbf{0} \\ * & * & * & -\eta_3 I \end{bmatrix} < 0, \tag{18}$$

in which

$$\begin{aligned} &\Sigma_{11} = sym\left\{\begin{bmatrix} A(rP_S + \bar{r}\overline{P_S}) + BT_4 & BT_3 \\ T_2 & T_1 \end{bmatrix}\right\}, \Sigma_{13} = \begin{bmatrix} (rP_S + \bar{r}\overline{P_S})^TA^T \\ \mathbf{0} \end{bmatrix}, \Sigma_{14} = \begin{bmatrix} T_4^TB & \mathbf{0} \\ T_3^TB & \mathbf{0} \end{bmatrix}, \Sigma_{33} = \varepsilon_1 G - I, \\ &\tilde{G} = \begin{bmatrix} G & \mathbf{0} \\ \mathbf{0} & \mathbf{0} \end{bmatrix}, \tilde{H} = \begin{bmatrix} H & \mathbf{0} \\ \mathbf{0} & \mathbf{0} \end{bmatrix}, \end{aligned} \tag{19}$$

where $\theta = (1 - \alpha)\pi/2$ then, the dynamic output feedback controller parameters of

$$A_C = T_1(rP_C + \bar{r}\overline{P_C})^{-1}, B_C = T_2(rP_S + \bar{r}\overline{P_S})^{-1}C^{\dagger}, C_C = T_3(rP_C + \bar{r}\overline{P_C})^{-1}, D_C = T_4(rP_S + \bar{r}\overline{P_S})^{-1}C^{\dagger}, \tag{20}$$

make the closed-loop system in (15) asymptotically stable.

**Proof**: It follows from Lemma 1 that the uncertain fractional-order closed-loop system (15) with $0 < \alpha < 1$ is asymptotically stable if there exists a positive definite matrix $P = P^*$, $P \in C^{(n+n_C)\times(n+n_C)}$ in the form of (17) in a way that

$$\begin{aligned}(rP+\bar{r}\bar{P})^T\tilde{A}_{Cl}^T+\tilde{A}_{Cl}(rP+\bar{r}\bar{P}) &= sym\left\{\begin{bmatrix} A(rP_S+\bar{r}\overline{P_S})+BD_CC(rP_S+\bar{r}\overline{P_S}) & BC_C(rP_C+\bar{r}\overline{P_C}) \\ B_CC(rP_S+\bar{r}\overline{P_S}) & A_C(rP_C+\bar{r}\overline{P_C})\end{bmatrix}\right\} \\ &+sym\left\{\begin{bmatrix} \Delta_A(rP_S+\bar{r}\overline{P_S}) & \mathbf{0} \\ \mathbf{0} & \mathbf{0}\end{bmatrix}\right\}+sym\left\{\begin{bmatrix} \Delta_I(A+\Delta_A)(rP_S+\bar{r}\overline{P_S}) & \mathbf{0} \\ \mathbf{0} & \mathbf{0}\end{bmatrix}\right\} \\ &+sym\left\{\begin{bmatrix} \Delta_IBD_CC(rP_S+\bar{r}\overline{P_S}) & \Delta_IBC_C(rP_C+\bar{r}\overline{P_C}) \\ \mathbf{0} & \mathbf{0}\end{bmatrix}\right\}<0,\end{aligned} \tag{21}$$

According to Lemma 3, the following three inequalities hold for any positive scalar constants $\eta_1$, $\eta_2$, and $\eta_3$

$$\begin{aligned} & sym\left\{\begin{bmatrix} \Delta_A(rP_S+\bar{r}\overline{P_S}) & \mathbf{0} \\ \mathbf{0} & \mathbf{0}\end{bmatrix}\right\} = sym\left\{\begin{bmatrix} \Delta_A & \mathbf{0} \\ \mathbf{0} & \mathbf{0}\end{bmatrix}\begin{bmatrix} (rP_S+\bar{r}\overline{P_S}) & \mathbf{0} \\ \mathbf{0} & \mathbf{0}\end{bmatrix}\right\} \le \eta_1\begin{bmatrix} \Delta_A\Delta_A^T & \mathbf{0} \\ \mathbf{0} & \mathbf{0}\end{bmatrix}+ \\ & \eta_1^{-1}\begin{bmatrix} (rP_S+\bar{r}\overline{P_S})^T(rP_S+\bar{r}\overline{P_S}) & \mathbf{0} \\ \mathbf{0} & \mathbf{0}\end{bmatrix}, \\ & sym\left\{\begin{bmatrix} \Delta_I(A+\Delta_A)(rP_S+\bar{r}\overline{P_S}) & \mathbf{0} \\ \mathbf{0} & \mathbf{0}\end{bmatrix}\right\} = sym\left\{\begin{bmatrix} \Delta_I & \mathbf{0} \\ \mathbf{0} & \mathbf{0}\end{bmatrix}\begin{bmatrix} (A+\Delta_A)(rP_S+\bar{r}\overline{P_S}) & \mathbf{0} \\ \mathbf{0} & \mathbf{0}\end{bmatrix}\right\} \le \eta_2\begin{bmatrix} \Delta_I\Delta_I^T & \mathbf{0} \\ \mathbf{0} & \mathbf{0}\end{bmatrix}+ \\ & \eta_2^{-1}\begin{bmatrix} P_S^T(A+\Delta_A)(A+\Delta_A)^T(rP_S+\bar{r}\overline{P_S}) & \mathbf{0} \\ \mathbf{0} & \mathbf{0}\end{bmatrix}, \\ & sym\left\{\begin{bmatrix} \Delta_IBD_CC(rP_S+\bar{r}\overline{P_S}) & \Delta_IBC_C(rP_C+\bar{r}\overline{P_C}) \\ \mathbf{0} & \mathbf{0}\end{bmatrix}\right\} = sym\left\{\begin{bmatrix} \Delta_I & \mathbf{0} \\ \mathbf{0} & \mathbf{0}\end{bmatrix}\begin{bmatrix} BD_CC(rP_C+\bar{r}\overline{P_C}) & BC_C(rP_C+\bar{r}\overline{P_C}) \\ \mathbf{0} & \mathbf{0}\end{bmatrix}\right\} \le \\ & \eta_3\begin{bmatrix} \Delta_I\Delta_I^T & \mathbf{0} \\ \mathbf{0} & \mathbf{0}\end{bmatrix}+\eta_3^{-1}\begin{bmatrix} BD_CC(rP_C+\bar{r}\overline{P_C}) & BC_C(rP_C+\bar{r}\overline{P_C}) \\ \mathbf{0} & \mathbf{0}\end{bmatrix}^T\begin{bmatrix} BD_CC(rP_C+\bar{r}\overline{P_C}) & BC_C(rP_C+\bar{r}\overline{P_C}) \\ \mathbf{0} & \mathbf{0}\end{bmatrix}. \end{aligned} \tag{22}$$

On the other hand, for any scalar $\varepsilon_1 > 0$ such that $I - \varepsilon_1 G > 0$, form (5), $I - \varepsilon_1\Delta_A\Delta_A^T > I - \varepsilon_1 G > 0$, which means that $(I - \varepsilon_1\Delta_A\Delta_A^T)^{-1} < (I - \varepsilon_1 G)^{-1}$. By using Lemma 4, we have

$$(A+\Delta_A)^T(A+\Delta_A) \le A^T(I-\varepsilon_1\Delta_A\Delta_A^T)^{-1}A+\varepsilon_1^{-1}I \le A^T(I-\varepsilon_1 G)^{-1}A+\varepsilon_1^{-1}I. \tag{23}$$

For simplicity let $\eta_2 = 1$ in (22), It easily follows from (5), and (21)-(23) that there exist scalar constant $\varepsilon_1 > 0$ in a way that

$$\begin{aligned}(rP+\bar{r}\bar{P})^T\tilde{A}_{Cl}^T+\tilde{A}_{Cl}(rP+\bar{r}\bar{P}) &\le sym\left\{\begin{bmatrix} A(rP_S+\bar{r}\overline{P_S})+BD_CC(rP_S+\bar{r}\overline{P_S}) & BC_C(rP_C+\bar{r}\overline{P_C}) \\ B_CC(rP_S+\bar{r}\overline{P_S}) & A_C(rP_C+\bar{r}\overline{P_C})\end{bmatrix}\right\}+\eta_1\begin{bmatrix} G & \mathbf{0} \\ \mathbf{0} & \mathbf{0}\end{bmatrix} \\ &+\eta_1^{-1}\begin{bmatrix} (rP_S+\bar{r}\overline{P_S})^T(rP_S+\bar{r}\overline{P_S}) & \mathbf{0} \\ \mathbf{0} & \mathbf{0}\end{bmatrix}+\begin{bmatrix} H & \mathbf{0} \\ \mathbf{0} & \mathbf{0}\end{bmatrix}+\begin{bmatrix} (rP_S+\bar{r}\overline{P_S})^TA^T(1-\varepsilon_1G)A(rP_S+\bar{r}\overline{P_S}) & \mathbf{0} \\ \mathbf{0} & \mathbf{0}\end{bmatrix}+\eta_3\begin{bmatrix} H & \mathbf{0} \\ \mathbf{0} & \mathbf{0}\end{bmatrix} \\ &+\eta_3^{-1}\begin{bmatrix} BD_CC(rP_C+\bar{r}\overline{P_C}) & BC_C(rP_C+\bar{r}\overline{P_C}) \\ \mathbf{0} & \mathbf{0}\end{bmatrix}^T\begin{bmatrix} BD_CC(rP_C+\bar{r}\overline{P_C}) & BC_C(rP_C+\bar{r}\overline{P_C}) \\ \mathbf{0} & \mathbf{0}\end{bmatrix}.\end{aligned} \tag{24}$$

By applying Schur complement to (24) one can obtain

$$\begin{bmatrix} \Sigma'_{11}+\eta_1\tilde{G}+(\eta_3+1)\tilde{H} & (rP+\bar{r}\bar{P})^T & \Sigma'_{13} & \Sigma'_{14} \\ * & -\eta_1 I & \mathbf{0} & \mathbf{0} \\ * & * & \Sigma'_{33} & \mathbf{0} \\ * & * & * & -\eta_3 I\end{bmatrix}<0, \tag{25}$$

in which

$$\begin{aligned}\Sigma'_{11} &= sym\left\{\begin{bmatrix} A(rP_S+\bar{r}\overline{P_S})+BD_CC(rP_S+\bar{r}\overline{P_S}) & BC_C(rP_C+\bar{r}\overline{P_C}) \\ B_CC(rP_S+\bar{r}\overline{P_S}) & A_C(rP_C+\bar{r}\overline{P_C})\end{bmatrix}\right\}, \quad \Sigma'_{13} = \begin{bmatrix} (rP_S+\bar{r}\overline{P_S})A^T & \mathbf{0} \\ \mathbf{0} & \mathbf{0}\end{bmatrix}, \\ \Sigma'_{14} &= \begin{bmatrix} (rP_S+\bar{r}\overline{P_S})C^TD_C^TB & \mathbf{0} \\ (rP_C+\bar{r}\overline{P_C})^TC_C^TB & \mathbf{0}\end{bmatrix}, \quad \Sigma'_{33} = \begin{bmatrix} \varepsilon_1 G - I & \mathbf{0} \\ \mathbf{0} & \mathbf{0}\end{bmatrix}\end{aligned} \tag{26}$$

Since the matrix inequality (26) is not linear owing to several multiplications of variables, the following change of variables

$$T_1 = A_C(rP_C+\bar{r}\overline{P_C}), \qquad T_2 = B_CC(rP_S+\bar{r}\overline{P_S}), \qquad T_3 = C_C(rP_C+\bar{r}\overline{P_C}), \qquad T_4 = D_CC(rP_S+\bar{r}\overline{P_S}) \tag{27}$$

solves the nonlinearity of the inequality (26) and also completes the proof. ■

And for $1 \le \alpha < 2$ we have the following theorem.

**Theorem 2**: Considering closed-loop system in (15) with $1 \le \alpha < 2$, if positive definite matrix $P = P^T$ in the form of (17) with $P_S \in \mathcal{R}^{n\times n}$ and $P_C \in \mathcal{R}^{n_c\times n_c}$ and real positive scalar constants $\eta_1$, $\varepsilon_1$, and $\eta_3$ alongside with matrices $T_i, i = 1, \dots, 4$ exist such that the following LMI constrain becomes feasible

$$\begin{bmatrix} \Sigma_{11}+\eta_1 I_2\otimes\tilde{G}+(\eta_3+1)I_2\otimes\tilde{H} & I_2\otimes P & \Sigma_{13} & \Sigma_{14} \\ * & -\eta_1 I & \mathbf{0} & \mathbf{0} \\ * & * & \Sigma_{33} & \mathbf{0} \\ * & * & * & -\eta_3 I \end{bmatrix} < 0, \tag{28}$$

in which

$$\Sigma_{11}=\begin{bmatrix} \sin\theta & \cos\theta \\ -\cos\theta & \sin\theta \end{bmatrix}\otimes\begin{bmatrix} AP_S+P_SA^T+BT_4+T_4^TB^T & BT_3+T_2^T \\ T_2+T_3^TB^T & T_1+T_1^T \end{bmatrix},\ \Sigma_{13}=I_2\otimes\begin{bmatrix} P_SA^T & \mathbf{0} \\ \mathbf{0} & \mathbf{0} \end{bmatrix},$$

$$\Sigma_{13}=I_2\otimes\begin{bmatrix} T_4^TB & \mathbf{0} \\ T_3^TB & \mathbf{0} \end{bmatrix},\ \Sigma_{33}=I_2\otimes\begin{bmatrix} \varepsilon_1 g-I & \mathbf{0} \\ \mathbf{0} & \mathbf{0} \end{bmatrix},\ \tilde{G}=\begin{bmatrix} G & \mathbf{0} \\ \mathbf{0} & \mathbf{0} \end{bmatrix},\ \tilde{H}=\begin{bmatrix} H & \mathbf{0} \\ \mathbf{0} & \mathbf{0} \end{bmatrix}, \tag{29}$$

where $\theta=\pi-\alpha\pi/2$ then, the dynamic output feedback controller parameters of

$$A_C=T_1{P_C}^{-1}, B_C=T_2{P_S}^{-1}C^{\uparrow}, C_C=T_3{P_C}^{-1}, D_C=T_4{P_S}^{-1}C^{\uparrow}, \tag{30}$$

make the closed-loop system in (15) asymptotically stable.

**Proof**: It follows from Lemma 2 that the uncertain fractional-order closed-loop system (15) with $1<\alpha\leq 2$ is asymptotically stable if there exists a positive definite matrix $P=P^T,\ \ P\in\mathcal{R}^{(n+n_C)\times(n+n_C)}$ in the form of (17) such that

$$\begin{bmatrix} \left(\tilde{A}_{Cl}P+P{\tilde{A}_{Cl}}^T\right)\sin\theta & \left(\tilde{A}_{Cl}P-P{\tilde{A}_{Cl}}^T\right)\cos\theta \\ \left(P{\tilde{A}_{Cl}}^T-\tilde{A}_{Cl}P\right)\cos\theta & \left(\tilde{A}_{Cl}P+P{\tilde{A}_{Cl}}^T\right)\sin\theta \end{bmatrix}=sym\left(\begin{bmatrix} \sin\theta & \cos\theta \\ -\cos\theta & \sin\theta \end{bmatrix}\otimes\begin{bmatrix} AP_S+BD_CCP_S & BC_CP_C \\ B_CCP_S & A_CP_C \end{bmatrix}\right)$$
$$+sym\left(\begin{bmatrix} \sin\theta & \cos\theta \\ -\cos\theta & \sin\theta \end{bmatrix}\otimes\begin{bmatrix} \Delta_AP_S & \mathbf{0} \\ \mathbf{0} & \mathbf{0} \end{bmatrix}\right)+sym\left(\begin{bmatrix} \sin\theta & \cos\theta \\ -\cos\theta & \sin\theta \end{bmatrix}\otimes\begin{bmatrix} \Delta_I(A+\Delta_A)P_S & \mathbf{0} \\ \mathbf{0} & \mathbf{0} \end{bmatrix}\right) \tag{31}$$
$$+sym\left(\begin{bmatrix} \sin\theta & \cos\theta \\ -\cos\theta & \sin\theta \end{bmatrix}\otimes\begin{bmatrix} \Delta_IBD_CCP_S & \Delta_IBC_CP_C \\ \mathbf{0} & \mathbf{0} \end{bmatrix}\right)<0,$$

According to Lemma 3, the following three inequalities hold for any positive scalar constants $\eta_1$, $\eta_2$, and $\eta_3$

$$sym\left(\begin{bmatrix} \sin\theta & \cos\theta \\ -\cos\theta & \sin\theta \end{bmatrix}\otimes\begin{bmatrix} \Delta AP_S & \mathbf{0} \\ \mathbf{0} & \mathbf{0} \end{bmatrix}\right)=sym\left(\left(\begin{bmatrix} \sin\theta & \cos\theta \\ -\cos\theta & \sin\theta \end{bmatrix}\otimes\begin{bmatrix} \Delta_A & \mathbf{0} \\ \mathbf{0} & \mathbf{0} \end{bmatrix}\right)\left(I_2\otimes\begin{bmatrix} P_S & \mathbf{0} \\ \mathbf{0} & \mathbf{0} \end{bmatrix}\right)\right)\leq\eta_1\left(I_2\otimes\begin{bmatrix} \Delta_A\Delta_A^T & \mathbf{0} \\ \mathbf{0} & \mathbf{0} \end{bmatrix}\right)+\eta_1^{-1}\left(I_2\otimes\begin{bmatrix} P_S^2 & \mathbf{0} \\ \mathbf{0} & \mathbf{0} \end{bmatrix}\right),$$

$$sym\left(\begin{bmatrix} \sin\theta & \cos\theta \\ -\cos\theta & \sin\theta \end{bmatrix}\otimes\begin{bmatrix} \Delta_I(A+\Delta_A)P_S & \mathbf{0} \\ \mathbf{0} & \mathbf{0} \end{bmatrix}\right)=sym\left(\left(\begin{bmatrix} \sin\theta & \cos\theta \\ -\cos\theta & \sin\theta \end{bmatrix}\otimes\begin{bmatrix} \Delta_I & \mathbf{0} \\ \mathbf{0} & \mathbf{0} \end{bmatrix}\right)\left(I_2\otimes\begin{bmatrix} (A+\Delta_A)P_S & \mathbf{0} \\ \mathbf{0} & \mathbf{0} \end{bmatrix}\right)\right)\leq$$
$$\eta_2\left(I_2\otimes\begin{bmatrix} \Delta_I\Delta_I^T & \mathbf{0} \\ \mathbf{0} & \mathbf{0} \end{bmatrix}\right)+\eta_2^{-1}\left(I_2\otimes\begin{bmatrix} P_S^T(A+\Delta_A)(A+\Delta_A)^TP_S & \mathbf{0} \\ \mathbf{0} & \mathbf{0} \end{bmatrix}\right), \tag{32}$$

$$sym\left(\begin{bmatrix} \sin\theta & \cos\theta \\ -\cos\theta & \sin\theta \end{bmatrix}\otimes\begin{bmatrix} \Delta_IBD_CCP_S & \Delta_IBC_CP_C \\ \mathbf{0} & \mathbf{0} \end{bmatrix}\right)=sym\left(\left(\begin{bmatrix} \sin\theta & \cos\theta \\ -\cos\theta & \sin\theta \end{bmatrix}\otimes\begin{bmatrix} \Delta_I & \mathbf{0} \\ \mathbf{0} & \mathbf{0} \end{bmatrix}\right)\left(I_2\otimes\begin{bmatrix} BD_CCP_C & BC_CP_C \\ \mathbf{0} & \mathbf{0} \end{bmatrix}\right)\right)\leq\eta_3\left(I_2\otimes\begin{bmatrix} \Delta_I\Delta_I^T & \mathbf{0} \\ \mathbf{0} & \mathbf{0} \end{bmatrix}\right)+\eta_3^{-1}\left(I_2\otimes\left(\begin{bmatrix} BD_CCP_C & BC_CP_C \\ \mathbf{0} & \mathbf{0} \end{bmatrix}^T\begin{bmatrix} BD_CCP_C & BC_CP_C \\ \mathbf{0} & \mathbf{0} \end{bmatrix}\right)\right).$$

Let $\eta_2=1$ in (32), It easily follows from (5), (23), (31), and (32) that there exist scalar constant $\varepsilon_1>0$ in a way that

$$\begin{bmatrix} \left(\tilde{A}_{Cl}P+P{\tilde{A}_{Cl}}^T\right)\sin\theta & \left(\tilde{A}_{Cl}P-P{\tilde{A}_{Cl}}^T\right)\cos\theta \\ \left(P{\tilde{A}_{Cl}}^T-\tilde{A}_{Cl}P\right)\cos\theta & \left(\tilde{A}_{Cl}P+P{\tilde{A}_{Cl}}^T\right)\sin\theta \end{bmatrix}\leq sym\left(\begin{bmatrix} \sin\theta & \cos\theta \\ -\cos\theta & \sin\theta \end{bmatrix}\otimes\begin{bmatrix} AP_S+BD_CCP_S & BC_CP_C \\ B_CCP_S & A_CP_C \end{bmatrix}\right)$$
$$+\eta_1\left(I_2\otimes\begin{bmatrix} G & \mathbf{0} \\ \mathbf{0} & \mathbf{0} \end{bmatrix}\right)+\eta_1^{-1}\left(I_2\otimes\begin{bmatrix} P_S^TP_S & \mathbf{0} \\ \mathbf{0} & \mathbf{0} \end{bmatrix}\right)+\left(I_2\otimes\begin{bmatrix} H & \mathbf{0} \\ \mathbf{0} & \mathbf{0} \end{bmatrix}\right)+\left(I_2\otimes\begin{bmatrix} P_S^TA^T(1-\varepsilon_1G)AP_S & \mathbf{0} \\ \mathbf{0} & \mathbf{0} \end{bmatrix}\right)+\eta_3\left(I_2\otimes\begin{bmatrix} H & \mathbf{0} \\ \mathbf{0} & \mathbf{0} \end{bmatrix}\right)+\eta_3^{-1}\left(I_2\otimes\begin{bmatrix} P_S^TD_C^TB^TBD_CP_S & P_S^TD_C^TB^TBD_CP_S \\ P_S^TC_C^TB^TBC_CP_S & P_S^TC_C^TB^TBC_CP_S \end{bmatrix}\right). \tag{33}$$

By applying Schur complement to (33) one can obtain

$$\begin{bmatrix} \Sigma'_{11}+\eta_1 G+(\eta_3+1)H & I_2\otimes P & \Sigma'_{13} & \Sigma'_{14} \\ * & -\eta_1 I & \mathbf{0} & \mathbf{0} \\ * & * & \Sigma'_{33} & \mathbf{0} \\ * & * & * & -\eta_3 I \end{bmatrix} < 0, \tag{34}$$

in which

$$\Sigma'_{11} = \begin{bmatrix} \sin\theta & \cos\theta \\ -\cos\theta & \sin\theta \end{bmatrix} \otimes \begin{bmatrix} AP_S + P_S A^T + BD_C CP_S + P_S C^T D_C^T B^T & BC_C P_C + P_S C^T B_C^T \\ B_C CP_S + T_3^T B^T & A_C P_C + P_C A_C \end{bmatrix}, \Sigma'_{13} = I_2 \otimes \begin{bmatrix} P_S A^T & \mathbf{0} \\ \mathbf{0} & \mathbf{0} \end{bmatrix},$$
$$\Sigma'_{13} = I_2 \otimes \begin{bmatrix} P_S C^T D_C^T B & \mathbf{0} \\ P_C^T C_C^T B & \mathbf{0} \end{bmatrix}, \Sigma'_{33} = I_2 \otimes \begin{bmatrix} \varepsilon_1 g - I & \mathbf{0} \\ \mathbf{0} & \mathbf{0} \end{bmatrix} \tag{35}$$

Nevertheless, the matrix inequality (35) is not linear owing to several multiplications of variables. Therefore, by linearizing change of variables as

$$T_1 = A_C P_C, T_2 = B_C CP_S, T_3 = C_C P_C, T_4 = D_C CP_S \tag{36}$$

the inequality in (34) turns into the one in (28) and it completes the proof. ■

**Remark 1**: The control scheme of the presented robust stabilization controller has been depicted in Fig. 1, in which the controller unknown parameters are calculated via proposed LMI constraints in Theorem 1 and 2, using system's uncertain and certain matrices. The outputs of the system are then fed to the controller and the control input is generated for uncertain FO-LTI system.

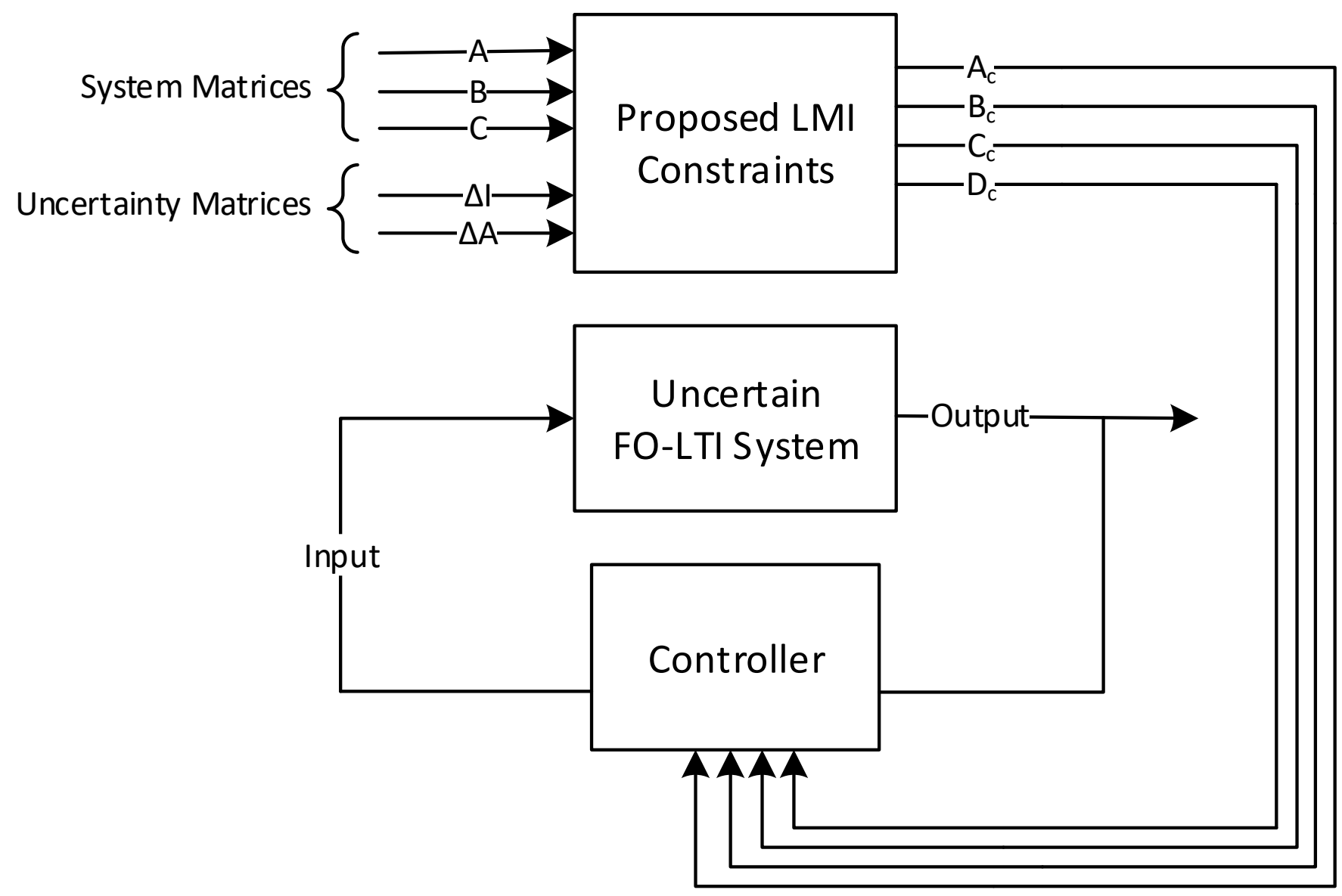

Fig. 1. The presented robust stabilization control scheme.

**Corollary 1**: However, In the following, two robust stabilizing theorems for uncertain FO-LTI system (12) are given in terms of LMIs for $0 < \alpha < 1$ and $1 \leq \alpha < 2$ cases, respectively.

**Theorem 1** and Theorem 2 are valid for robust stabilization of uncertain FO-LTI systems of form (2), the proposed methods can also be easily used for the case of certain systems by solving the LMI constraints $\Sigma_{11} < 0$ in these theorems.

**Proof:** The proof is simple by assuming $\Delta_I = \Delta_A = \mathbf{0}$ in the proof procedure of Theorem 1 and Theorem 2.

**Remark 2**: In Theorem 1 and Theorem 2 by checking the feasibility of the LMI conditions, unknown parameters of the stabilizing dynamic output feedback controllers can be achieved.

**Remark 3**: According to (14) the proposed dynamic output feedback controller contains the $\boldsymbol{D_C}$ matrix which is "the direct feedthrough parameter", that can be obtained from (27) and (36) in the proof procedure of the Theorem 1 and Theorem 2, respectively.

## 4. Numerical examples

This section provides some numerical examples that illustrate the effectiveness of the proposed methods in this paper. In this paper, we use YALMIP parser (Löfberg, 2004) and SeDuMi (Sturm, 1999) solver in Matlab tool (Higham and Higham, 2005) to check the feasibility of the proposed constraints and obtain the controller parameters. Grünwald-Letnikov derivative is used to discretize the fractional-order derivatives in simulating the numerical examples.

### 4.1. Example 1 *for the* $0 < \alpha < 1$ *case*

The following fractional-order linear system with uncertain parameters $\Delta_I$ and $\Delta_A$ is considered (Chen et al., 2015)

$$D^\alpha x(t) = (I + \Delta I)[(A + \Delta A)x(t) + Bu(t)]. \tag{37}$$

Consider the dynamic output feedback stabilization problem of the uncertain fractional-order system (12) with $\alpha = 0.65$, where

$$A = \begin{bmatrix} -2 & 0 & -1 \\ 0 & 3 & 0 \\ -1 & -1 & 4 \end{bmatrix}, \Delta_I = \begin{bmatrix} 0 & r_1 & 0 \\ 0 & 0 & 0 \\ 0 & 0 & r_2 \end{bmatrix}, \Delta_A = \begin{bmatrix} 0 & 0 & 0 \\ 0 & s_1 & s_2 \\ 0 & s_1 & s_2 \end{bmatrix}, B = \begin{bmatrix} 1 & 0 & 0 \\ 0 & 1 & 0 \\ 0 & 0 & 1 \end{bmatrix}, C = \begin{bmatrix} 1 & 0 & 1 \end{bmatrix}, \tag{38}$$

$|r_i| \leq 0.3, \qquad |s_i| \leq 0.3, \qquad i = 1,2.$

by some calculation, it can be obtained that

$$H = \begin{bmatrix} 0.18 & 0 & 0 \\ 0 & 0 & 0 \\ 0 & 0 & 0.18 \end{bmatrix}, \qquad G = \begin{bmatrix} 0 & 0 & 0 \\ 0 & 0.36 & 0.36 \\ 0 & 0.36 & 0.36 \end{bmatrix}. \tag{39}$$

The eigenvalues of $\tilde{A}$, $A_{Cl}$, and stability boundaries $\pm\alpha\pi/2$ are depicted in Fig. 2. According to Lemma 1 and Fig. 2, the system (12) with parameters in (39) is unstable due to some eigenvalues of $\tilde{A}$ which are located on the right side of the stability boundaries. Yet, according to In the following, two robust stabilizing theorems for uncertain FO-LTI system (12) are given in terms of LMIs for $0 < \alpha < 1$ and $1 \leq \alpha < 2$ cases, respectively.

**Theorem 1**, it can be concluded that this uncertain fractional-order system is asymptotically stabilizable using the obtained dynamic output feedback controllers of arbitrary orders in the form of (14), listed in Table 1. The eigenvalues of $A_{Cl}$ are placed in the stability region which is also evident in Fig. 2.

The state trajectories of the resulted uncertain closed-loop FO-LTI system of the form (15), through obtained controllers with $n_C = 1$ and the static controller introduced in (Chen et al., 2015) are plotted in Fig. 3, where all states asymptotically converge to zero. It can be seen that the obtained dynamic output feedback controllers, even with a lower order of $n_C = 1$, have more efficient stabilizing actions compared to the static one. The settling time of the closed-loop system via proposed controller is very small compared to the static controller proposed in (Chen et al., 2015) as it is obvious from Fig. 3.

Table 1. Controller parameters obtained by Theorem 1 for Example 1.

| $n_C$ | $A_C$ | $B_C$ | $C_C$ | $D_C$ |
|---|---|---|---|---|
| **0** | 0 | 0 | 0 | $\begin{bmatrix} -1.813 & 0.417 & 4.216 \\ 0.680 & -12.002 & -4.945 \\ 4.721 & -4.005 & -30.816 \end{bmatrix}$ |
| **1** | $-1.2453$ | $\begin{bmatrix} -0.0246 \\ 0.2587 \\ 0.1423 \end{bmatrix} \times 10^{-3}$ | $\begin{bmatrix} 0.003 \\ -0.536 \\ -0.296 \end{bmatrix} \times 10^{-4}$ | $\begin{bmatrix} -1.825 & 0.558 & 4.602 \\ 0.917 & -12.830 & -6.526 \\ 5.224 & -5.342 & -34.139 \end{bmatrix}$ |
| **2** | $\begin{bmatrix} -1.2197 & 0 \\ 0 & -1.2197 \end{bmatrix}$ | $\begin{bmatrix} -0.0759 & -0.0078 \\ 0.1270 & 0.1375 \\ 0.1503 & 0.0145 \end{bmatrix} \times 10^{-3}$ | $\begin{bmatrix} 0.314 & 0.064 \\ -0.230 & -0.188 \\ 0.122 & 0.151 \end{bmatrix} \times 10^{-4}$ | $\begin{bmatrix} -1.847 & 0.698 & 4.952 \\ 1.139 & -13.598 & -8.014 \\ 5.680 & -6.644 & -37.184 \end{bmatrix}$ |

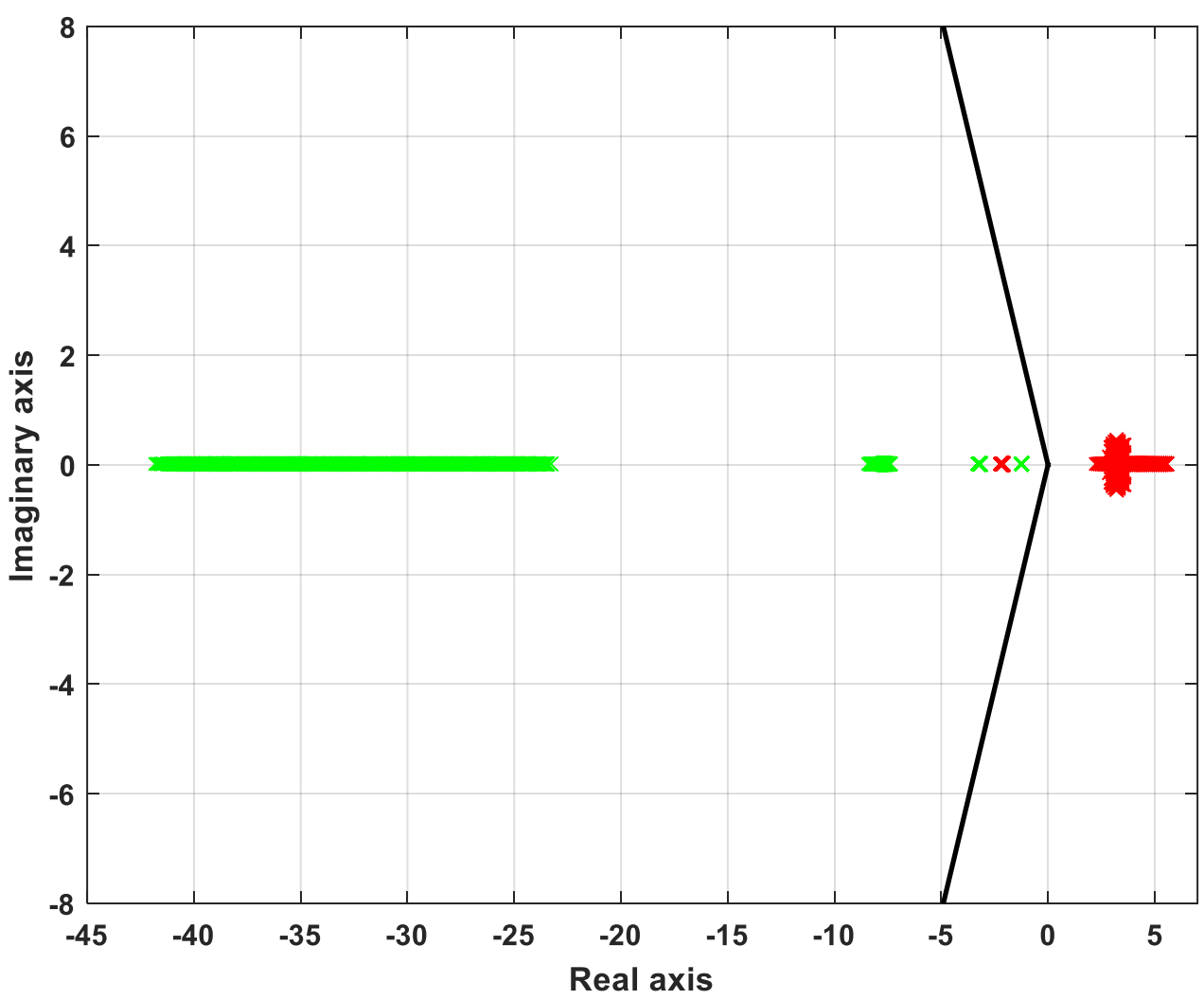

Fig. 2. The location of eigenvalues of the uncertain open-loop system (red) and closed-loop system via obtained output feedback controller with $n_C = 1$ (green) in Example 1.

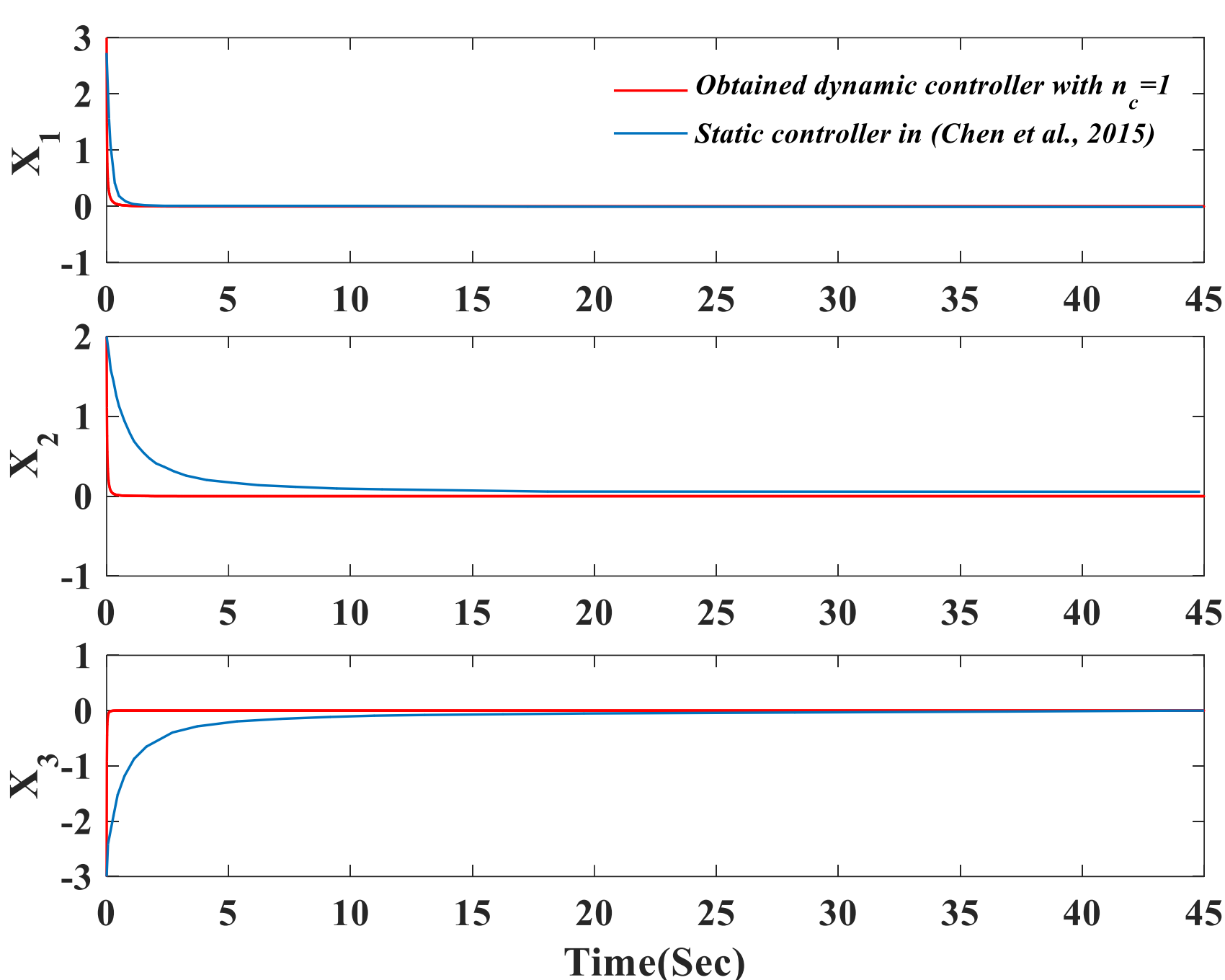

Fig. 3. The time response of the closed-loop system in Example 1 via obtained output feedback controllers with $n_C = 1$ and the static controller in (Chen et al., 2015).

4.2. Example *2 for the* $1 \leq \alpha < 2$ *case*

Dynamic output feedback stabilization problem of the uncertain fractional-order system (12) is considered with the parameters in (39) and $\alpha = 1.25$ (Chen et al., 2015).

The eigenvalues of $\tilde{A}$, $A_{Cl}$, and stability boundaries $\pm\alpha\pi/2$ are plotted in Fig. 4. According to Lemma 2 and Fig. 4, the system (12) with parameters in (39) and $\alpha = 1.25$ is unstable because of some eigenvalues of $\tilde{A}$ which are located on the right side of the stability boundaries. However, according to Theorem 2, it can be deduced that this uncertain fractional-order system is asymptotically stabilizable employing the obtained dynamic output feedback controllers of the form (14) and with arbitrary orders, tabulated in Table 2. All of the eigenvalues of $A_{Cl}$ are in the stability region which is also obvious in Fig. 4.

The state trajectories of the resulted uncertain closed-loop FO-LTI system of form the (15), through obtained controllers with $n_C = 1$ and the static controller introduced in (Chen et al., 2015) are plotted in Fig. 5, where all the states asymptotically converge to zero. It can be seen that the obtained dynamic output feedback controllers, even with a low order of $n_C = 1$, have more efficient stabilizing actions compared to the static one. The settling time of the closed-loop system via proposed controller is very small compared to the static controller proposed in (Chen et al., 2015) as it is obvious from Fig. 5.

Table 2. Controller parameters obtained by Theorem 2 for Example 2.

| $n_C$ | $A_C$ | $B_C$ | $C_C$ | $D_C$ |
|---|---|---|---|---|
| **0** | 0 | 0 | 0 | $\begin{bmatrix} 1.068 & -3.248 & 8.138 \\ -17.141 & 2.089 & -14.099 \\ -10.091 & 9.500 & -62.273 \end{bmatrix}$ |
| **1** | $-16.067$ | $\begin{bmatrix} 0.133 \\ 0.075 \\ -0.189 \end{bmatrix} \times 10^{-3}$ | $\begin{bmatrix} 0.289 \\ -0.298 \\ -0.296 \end{bmatrix} \times 10^{-3}$ | $\begin{bmatrix} 0.836 & -15.363 & 7.319 \\ -19.909 & 1.619 & -12.812 \\ -8.994 & 8.571 & -73.459 \end{bmatrix}$ |
| **2** | $\begin{bmatrix} -1.2197 & 0 \\ 0 & -1.2197 \end{bmatrix}$ | $\begin{bmatrix} -0.0759 & -0.0078 \\ 0.1270 & 0.1375 \\ 0.1503 & 0.0145 \end{bmatrix} \times 10^{-3}$ | $\begin{bmatrix} 0.314 & 0.064 \\ -0.230 & -0.188 \\ 0.122 & 0.151 \end{bmatrix} \times 10^{-4}$ | $\begin{bmatrix} -1.847 & 0.698 & 4.952 \\ 1.139 & -13.598 & -8.014 \\ 5.680 & -6.644 & -37.184 \end{bmatrix}$ |

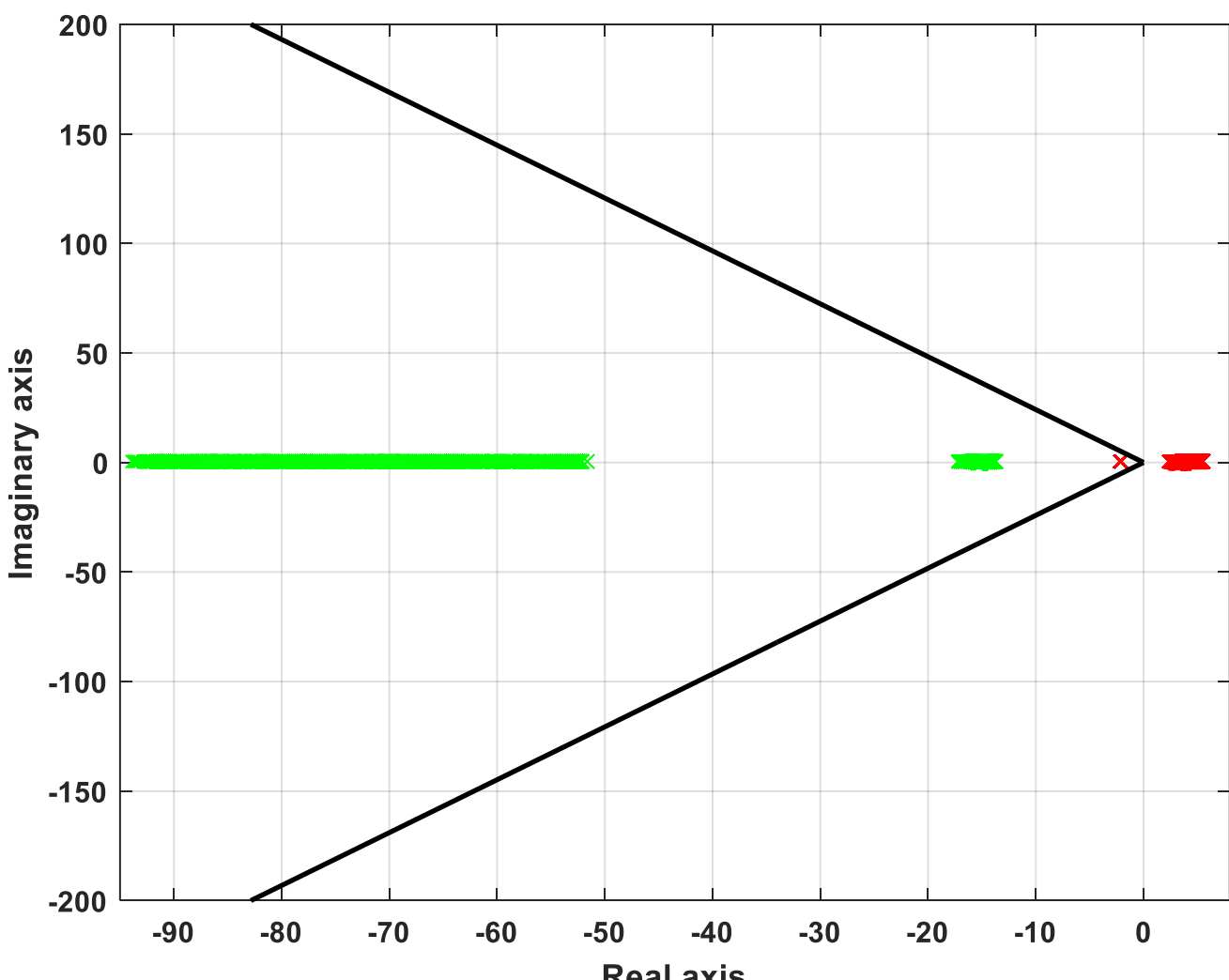

Fig. 4. The location of eigenvalues of the uncertain open-loop system (red) and closed-loop system via obtained output feedback controller with $n_C = 1$ (green) in Example 2.

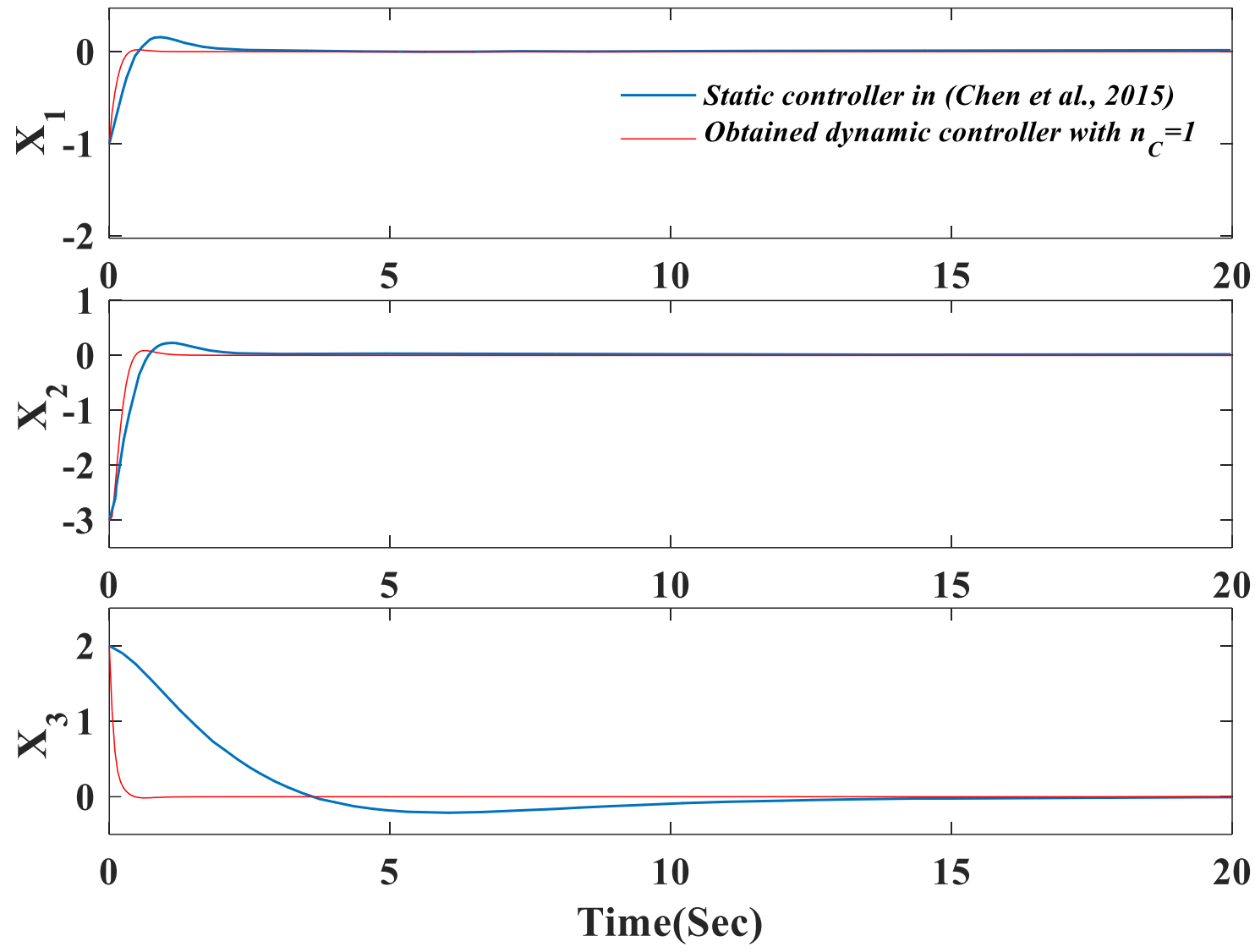

Fig. 5. The time response of the closed-loop system in Example 2 via obtained output feedback controllers with $n_C = 1$ and the static controller in (Chen et al., 2015).

4.3. Example *3 for the* $1 \le \alpha < 2$ *case*

Dynamic output feedback stabilization problem of the uncertain fractional-order system (12) is considered with the following parameters and $\alpha = 1.5$ (Xing and Lu, 2009)

$$A = \begin{bmatrix} 0 & 1 & 0 \\ 0 & 0 & 1 \\ 1 & -4 & 4 \end{bmatrix}, \Delta_I = \begin{bmatrix} 0 & 0 & 0 \\ 0 & 0 & 0 \\ 0 & 0 & r_1 \end{bmatrix}, \Delta_A = \begin{bmatrix} 0 & 0 & 0 \\ 0 & 0 & 0 \\ s_1 & s_2 & s_3 \end{bmatrix}, B = \begin{bmatrix} 0 \\ 0 \\ 1 \end{bmatrix}, C = \begin{bmatrix} 1 & 0 & 0 \\ 0 & 1 & 0 \\ 0 & 0 & 1 \end{bmatrix}, \tag{40}$$

$|r_1| \le 0{:}0157$, $|s_1| \le 0.01367$, $|s_2| \le 0{:}01933$ , $|s_3| \le 0{:}07135$.

by some calculations, it can be obtained that

$$H = \begin{bmatrix} 0 & 0 & 0 \\ 0 & 0 & 0 \\ 0 & 0 & 0.2465 \end{bmatrix} \times 10^{-3}, \quad G = \begin{bmatrix} 0 & 0 & 0 \\ 0 & 0 & 0 \\ 0 & 0 & 0.5651 \end{bmatrix} \times 10^{-2}. \tag{41}$$

The eigenvalues of $\tilde{A}$, $A_{Cl}$, and stability boundaries $\pm\alpha\pi/2$ are plotted in Fig. 6. According to Lemma 2 and Fig. 6, the system (12) with parameters in (40) and $\alpha = 1.5$ is unstable due to some eigenvalues of $\tilde{A}$ which are located on the right side of the stability boundaries. However, according to Theorem 2, it can be deduced that this uncertain fractional-order system is asymptotically stabilizable employing the obtained dynamic output feedback controllers of the form (14) and with arbitrary orders, tabulated in Table 3. All of the eigenvalues of $A_{Cl}$ are in the stability region which is also obvious in Fig. 6.

The state trajectories of the resulted uncertain closed-loop FO-LTI system of the form (15), through obtained controllers with $n_C = 1$ and the static controller introduced in (Xing and Lu, 2009) are plotted in Fig. 7, where all the states asymptotically converge to zero. It can be seen that the obtained dynamic output feedback controllers, even with a low order of $n_C = 1$, have more efficient stabilizing actions compared to the static one. The settling time of the closed-loop system via proposed controller is very small compared to the static controller proposed in (Xing and Lu, 2009) as it is obvious from Fig. 7.

Table 3. Controller parameters obtained by Theorem 2 for Example 3.

| $n_C$ | $A_C$ | $B_C$ | $C_C$ | $D_C$ |
|---|---|---|---|---|
| **0** | 0 | 0 | 0 | $[-124.818 \quad -37.795 \quad -63.821]$ |
| **1** | $-16.220$ | $\begin{bmatrix} 0.117 \\ 0.034 \\ 0.048 \end{bmatrix} \times 10^{-3}$ | $0.502 \times 10^{-4}$ | $[-83.379 \quad -25.279 \quad -43.617]$ |
| **2** | $\begin{bmatrix} -12.470 & 0 \\ 0 & -12.470 \end{bmatrix}$ | $\begin{bmatrix} -0.077 & 0.125 \\ -0.023 & 0.035 \\ -0.032 & 0.054 \end{bmatrix} \times 10^{-3}$ | $\begin{bmatrix} -0.240 \\ 0.096 \end{bmatrix} \times 10^{-4}$ | $[-76.775 \quad -23.341 \quad -40.940]$ |

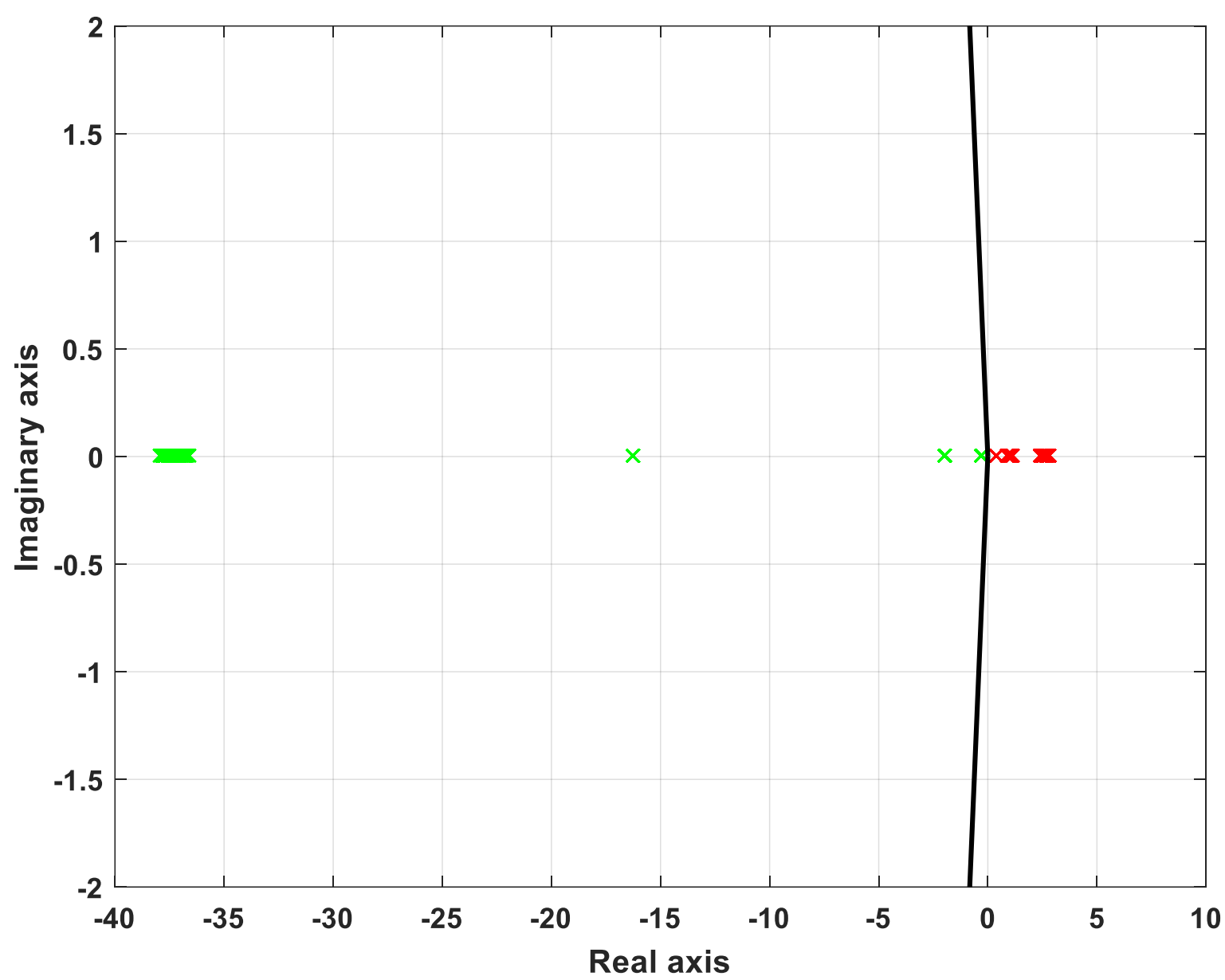

Fig. 6. The location of eigenvalues of the uncertain open-loop system (red) and closed-loop system via obtained output feedback controller with $n_C = 1$ (green) in Example 3.

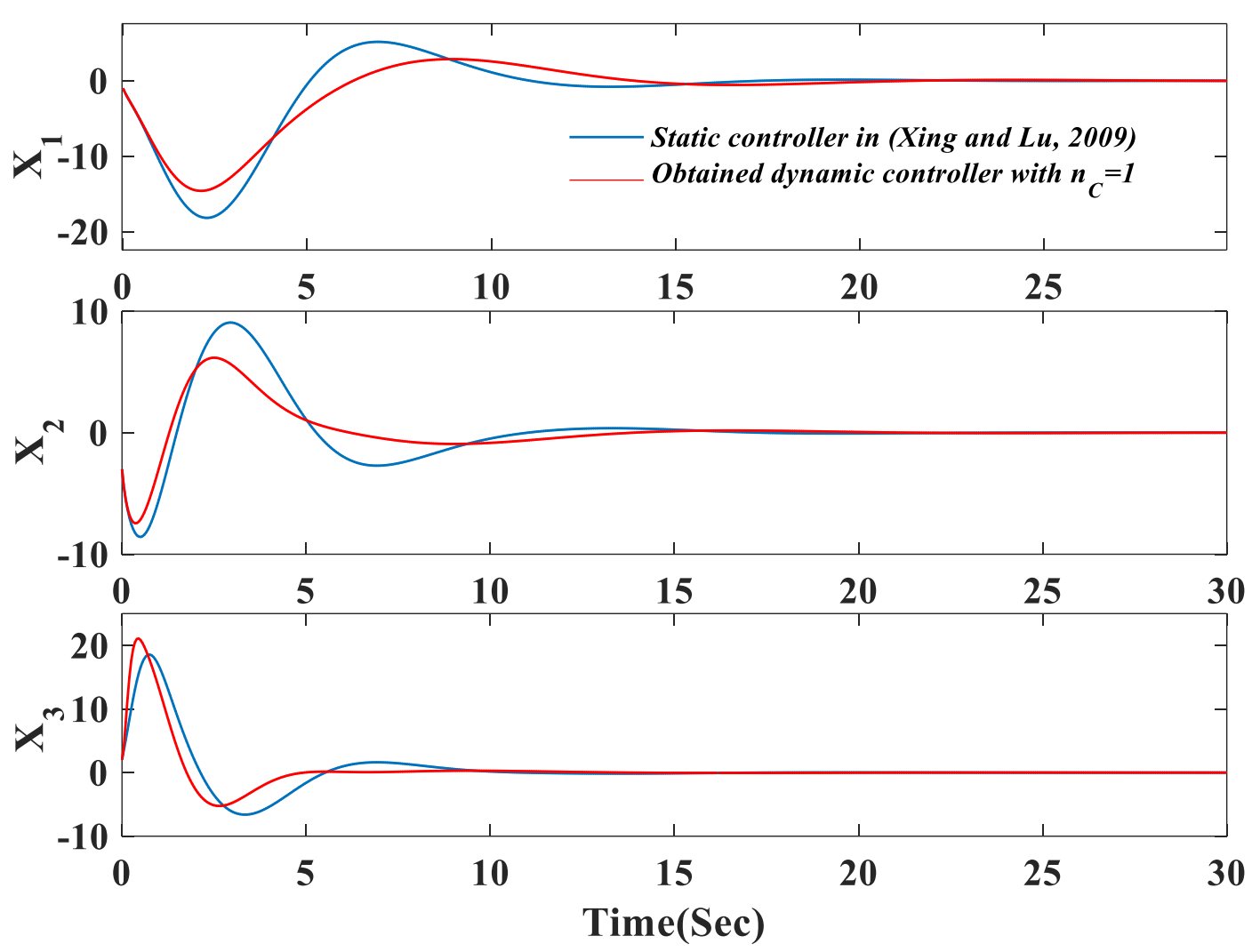

Fig. 7. The time response of the closed-loop system in Example 3 via obtained output feedback controllers with $n_C = 1$ and the static controller in (Xing and Lu, 2009).

## 5. Conclusion

In this paper the problem of robust dynamic output stabilization of uncertain FO-LTI systems with the fractional order $0 < \alpha < 2$, is solved in terms of LMIs. Sufficient conditions are obtained for designing a stabilizing controller with a predetermined order, which can be chosen to be as low as possible for simpler implementation. Indeed, by using the proposed method, one can benefit from dynamic output feedback controller advantages with lower orders than the system order. The LMI-based procedures of developing robust stabilizing control are preserved in spite of the complexity of assuming the most complete model of the linear controller, with a direct feedthrough parameter. Designing a dynamic robust controller, in order to take advantages of such controllers, leads to more unknown parameters in comparison with a static controller and makes controller design procedure more difficult due to more complex constraints which must be solved. In this paper, using proper lemmas and theorems, LMI techniques, and suitable solvers and parsers the difficulty of designing such controllers has been overcome.

Eventually, some numerical examples have shown the correctness and effectiveness of our results.